\documentclass[11pt,oneside]{amsart}

\usepackage{amsfonts,amsmath,amssymb,amsthm}
\usepackage{geometry}
\usepackage{graphicx,xcolor}
\usepackage[hidelinks]{hyperref}
\usepackage{dutchcal}

\usepackage[alphabetic,abbrev]{amsrefs} 
\renewcommand{\MR}[1]{}

\usepackage{zref-clever} \newcommand\Cref[1]{\zcref[S]{#1}}

\newcommand\newtheoremx[2]{%
\AddToHook{env/#1/begin}{
\zcsetup{countertype={theorem=#1}}}
\newtheorem{#1}[theorem]{#2}
}

\newtheorem{theorem}{Theorem}[section]
\newtheoremx{claim}{Claim}
\newtheoremx{proposition}{Proposition}
\newtheoremx{corollary}{Corollary}
\newtheoremx{lemma}{Lemma}
\newtheoremx{conjecture}{Conjecture}

\newtheorem{introtheorem}{Theorem}

\zcsetup{countertype={introtheorem=theorem}}

\theoremstyle{definition}
\newtheoremx{definition}{Definition}
\newtheoremx{notation}{Notation}
\newtheoremx{assumption}{Assumption}

\theoremstyle{remark}
\newtheoremx{example}{Example}
\newtheoremx{remark}{Remark}
\newtheoremx{question}{Question}

\numberwithin{equation}{section}

\newcommand\op{^{\mathrm{op}}}
\newcommand\cA{\mathcal A}
\newcommand\cB{\mathcal B}
\newcommand\cC{\mathcal C}
\newcommand\cD{\mathcal D}
\newcommand\cF{\mathcal F}
\newcommand\cP{\mathcal P}
\newcommand\cT{\mathcal T}
\newcommand\Ver{\mathrm{Ver}}
\newcommand\SL{\mathrm{SL}}
\newcommand\GL{\mathrm{GL}}
\newcommand\Vect{\mathrm{Vect}}
\newcommand\SVect{\mathrm{SVect}}
\newcommand\Dim{\operatorname{Dim}}
\newcommand\Rep{\operatorname{Rep}}
\newcommand\End{\operatorname{End}}

\newcommand\Tilt{\operatorname{Tilt}}
\newcommand\im{\operatorname{im}}
\newcommand\Sym{\operatorname{Sym}}
\newcommand\modfg{\operatorname{mod}}
\newcommand\gr{\operatorname{gr}}
\renewcommand\o\otimes
\renewcommand\:\colon
\newcommand\ol{\overline}
\newcommand\ti\widetilde
\newcommand\one{\mathbf{1}}
\newcommand\kk{\Bbbk}
\newcommand\id{\operatorname{id}}
\newcommand\Id{\operatorname{Id}}
\newcommand\Hom{\operatorname{Hom}}
\newcommand\sgd{\operatorname{sgd}}
\newcommand\len{\operatorname{len}}
\newcommand\FPdim{\operatorname{FPdim}}
\newcommand\sh{^{\mathrm{sh}}}
\newcommand\ZZ{\mathbb Z}
\newcommand\Ffrak{\mathfrak F}

\title{Maximal nilpotence for higher Verlinde categories}

\author{Kevin Coulembier}
\email{kevin.coulembier@sydney.edu.au}

\author{Johannes Flake}
\email{flake@math.uni-bonn.de}

\author{Joseph Newton}
\email{j.newton@maths.usyd.edu.au}

\author{Victor Ostrik}
\email{vostrik@uoregon.edu}

\begin{document}

\begin{abstract} We show that all higher Verlinde categories are maximally nilpotent, as conjectured by Coulembier (2023). By previous results, it suffices to show that non-trivial simple objects have only finitely many non-zero symmetric powers in odd characteristic. We compute the maximal degrees of these non-zero symmetric and exterior powers explicitly, hence confirming their finiteness. We discuss some of the far-reaching consequences.
\end{abstract}

\maketitle

\section{Introduction}

In \cite{Cou-commalg}, the notion of \emph{maximally nilpotent} (symmetric) tensor category was introduced in order to develop a tensor categorical version of commutative algebra in positive characteristic. Since then, the notion has proven extremely useful not only for establishing commutative algebra and algebraic geometry relative to a suitable general base category, but also for the classification of tensor categories in positive characteristic (see \Cref{sec:incompressible} below). Of particular interest in this context are the higher Verlinde categories $\Ver_{p^r}$, for $p$ a prime and $r\ge1$, as introduced in \cites{BEO,Cou-abelian}. They were shown to be analogs, in a certain sense, of the categories of vector spaces or super vector spaces, which are the underlying tensor categories used for commutative algebra and algebraic geometry in characteristic $0$. In particular, it was established recently \cite{CF-beo} that every finite tensor category over an algebraically closed field of positive characteristic fibers over a higher Verlinde category, as conjectured even for tensor categories of moderate growth in \cite{BEO}*{Conjecture~1.4}.

However, maximal nilpotence for the higher Verlinde categories remained a conjecture \cite{Cou-commalg}*{Conjecture~3.2.3}, up to the case $p=2$ confirmed in \cite{CEO-incompressible}. We show that the conjecture holds in general.

\begin{introtheorem}[\Cref{cor:MN}] \label{thm:A} The higher Verlinde categories $\Ver_{p^r}$, for all primes $p$ and all $1\le r\le \infty$, are maximally nilpotent.
\end{introtheorem}
    
By results in \cite{CEO-incompressible}, the conjecture was already reduced to showing that in odd characteristic all simple objects in the higher Verlinde categories other than the tensor unit have only finitely many non-zero symmetric powers. For comparison, the same statement for the category of super vector spaces follows from noting that the symmetric square of the odd one-dimensional space vanishes. We verify this for all $\Ver_{p^r}$ by deriving, in \Cref{thm:formulas} and \Cref{thm:formula-p2}, explicit formulas for the highest non-vanishing symmetric and exterior degrees, $m(Y)$ and $n(Y)$, for all simple objects $Y$. In particular, it follows that, if $Y$ is not isomorphic to the tensor unit or the odd tensor unit, then
$$
m(Y) + n(Y) \le p^r ,
$$
see \Cref{cor:bounds}, which confirms an upper bound suggested in \cite{CEN}, see \Cref{rem:correction}.

Our main technique centers around functors
$$
G_0\: \cD_{r,a}\to\Ver_{p^{r+1}} ,\quad
G\: [\cP_r\op,\cD_{r,a}] \to \Ver_{p^{r+1}} ,
$$
introduced in \Cref{sec:transfer}, where $\cD_{r,a}$ is a certain full additive and idempotent complete monoidal subcategory of $\Rep\GL_{a+1}$ and the source of $G$ is a non-abelian version of a Deligne tensor product of $\cD_{r,a}$ and $\Ver_{p^r}$. The functors $G$ and $G_0$ allow us to compute the relevant symmetric and exterior powers.

We discuss the consequences of \Cref{thm:A} in \Cref{sec:MN}.

We also consider the related notion of \emph{alternating powers} in \Cref{sec:alternating}. In particular, we show that for any non-zero object in a symmetric tensor category of moderate growth over any algebraically closed field of positive characteristic, there is a finite alternating power that is an invertible object with trivial self-braiding, while all higher alternating powers vanish.

Finally, we record some implications on (integrality of) $p$-adic dimensions and Hilbert series in \Cref{sec:dim}.

\subsection*{Use of AI} {GPT-6} Astra was used to obtain first versions of the proofs in this paper. We checked and simplified the proofs, and wrote this article by hand, before checking the document using the same LLM. 

\subsection*{Acknowledgments} We thank Pavel Etingof for crucial suggestions and questions throughout the project. K.C.~was supported by ARC grants FT220100125 and DP250100762. J.F.~was supported by the Hausdorff Center for Mathematics Bonn. J.N.~was supported by an Australian Government Research Training Program Scholarship.

\section{Background}

Let $\kk$ be an algebraically closed field of positive characteristic $p$.

\subsection{Tensor categories} In this paper, tensor categories are \emph{symmetric tensor categories} in the sense of \cite{EGNO}. In particular, they are hom-finite finite length abelian $\kk$-linear categories that are also symmetric rigid monoidal. This is consistent with the convention in \cite{CEN}. A \emph{tensor functor} is an exact linear symmetric monoidal functor between tensor categories. It is then automatically faithful. We will use \emph{Deligne tensor products} $\cC\boxtimes\cD$ of tensor categories $\cC,\cD$, and refer to \cite{CF-frob}*{Section~3} for a detailed discussion.

For a linear algebraic group $G$, we use the notation $\Rep G$ for the category of finite-dimensional rational $G$-modules over $\kk$, and this has the structure of a symmetric tensor category.

\subsection{Symmetric and exterior powers} Assume $\kk$ is of characteristic $p>0$. We recall some definitions and facts about symmetric and exterior powers from \cite{CEN}*{Section~2.1}. 

For any object $X$ in a tensor category in the above sense, we set
$$
\Lambda^2 X := \im(\id-c_{X,X}) \subset X^{\o2} ,\quad
\Gamma^2 X := \ker(\id-c_{X,X}) \subset X^{\o2} .
$$
For any $d\ge0$, the \emph{symmetric} and \emph{exterior powers} of $X$ are the quotient objects of $X^{\o d}$ defined by $\Sym^0 X:=\wedge^0 X:=\one$, $\Sym^1 X:=\wedge^1 X:=X$, and
$$
\Sym^d X := X^{\o d} \Big/ \sum_{i=0}^{d-2} X^{\o i} \o \Lambda^2 X \o X^{\o d-2-i} ,\quad
\wedge^d X := X^{\o d} \Big/ \sum_{i=0}^{d-2} X^{\o i} \o \Gamma^2 X \o X^{\o d-2-i} .
$$
The \emph{divided} and \emph{dual exterior powers} of $X$ are the subobjects of $X^{\o d}$ defined as
$$
\Gamma^d X:=(\Sym^d(X^*))^* ,\quad
\Lambda^d X:=(\wedge^d(X^*))^* .
$$
For $d=2$, this agrees with the initial definitions above. For $X\in\Rep G$, we have $\Lambda^d X\cong\wedge^d X$ in every degree.

For objects $X,Y$ in a tensor category and $d\ge0$, we have
\begin{equation}
\Sym^d (X\oplus Y)\cong \bigoplus_{0\le i\le d} \Sym^i X \o \Sym^{d-i} Y  , \quad
\Lambda^d (X\oplus Y)\cong \bigoplus_{0\le i\le d} \Lambda^i X \o \Lambda^{d-i} Y ,
\label{eq:direct-sums}
\end{equation}
by \cite{CEO-incompressible}*{Equation~3.10} and \cite{CEN}*{Theorem~2.2.5}, using duality for $\Lambda$. The second formula also holds with $\Lambda$ replaced by $\wedge$.

The first few symmetric and exterior powers are ``easier'' than the remaining ones, in the following sense.

\begin{lemma} For $0\le d\le p-1$, the objects $\Sym^d X$ and $\Lambda^d X$ are direct summands of $X^{\o d}$ given by idempotent endomorphisms in the image of the algebra map $\kk S_d\to\End(X^{\o d})$ arising from the braiding. Namely, these idempotents are the images of the symmetrizer $\frac{1}{d!}\sum_{\sigma\in S_d}\sigma$ and the antisymmetrizer $\frac{1}{d!}\sum_{\sigma\in S_d}(-1)^\sigma \sigma$, respectively. In these degrees, $\Lambda^d X\cong\wedge^d X$.
\end{lemma}

In particular, such symmetric and exterior powers are preserved by all linear symmetric monoidal functors, without requiring exactness.

\medskip

The \emph{symmetric degree} $m(X)$ and the \emph{exterior degree} $n(X)$ are defined as the maximal degrees of non-vanishing symmetric and dual exterior powers of $X$, respectively, with value $\infty$ if the degrees are unbounded. For instance, in any tensor category and for any $d\ge1$, $m(\one^d)=\infty$ and $n(\one^d)=d$.

Since $(\wedge^d X)^*\cong\Lambda^d(X^*)$, the maximal degree of a non-zero exterior power $\wedge^d X$ is $n(X^*)$.

The following criterion will be helpful for us to determine symmetric and exterior degrees.

\begin{lemma}[\cite{CEN}*{Proposition~2.3.2}] \label{lem:invertible} Let $X$ be non-invertible and $d\geq1$. Then $\Sym^d X$ is invertible if and only if $d=m(X)<\infty$, and $\Lambda^d X$ is invertible if and only if $d=n(X)<\infty$.
\end{lemma}

\subsection{Higher Verlinde categories} \label{bg:ver}

We recall the construction and some structural properties of the higher Verlinde categories from \cites{BEO,Cou-abelian}. Let $\Tilt\SL_2\subset\Rep\SL_2$ denote the category of tilting representations. It is an additive idempotent complete rigid symmetric monoidal subcategory of the tensor category $\Rep\SL_2$. The indecomposable objects $T_i$ in $\Tilt\SL_2$ can be labeled by their highest weight $i\ge0$. In particular, $T_0\cong\one$, $T_1$ is the natural 2-dimensional $\SL_2$-representation, and for all $i\ge0$, $T_i$ is the unique indecomposable summand of $T_1^{\o i}$ that is not a summand of $T_1^{\o j}$ for any $j<i$.

For any $r\ge0$, set $St_r:=T_{p^r-1}\in\Tilt\SL_2$, the $r$-th Steinberg representation of $\SL_2$. The non-zero tensor ideals / thick tensor ideals in $\Tilt\SL_2$ are exactly those generated by the objects $St_r$, for $r\ge0$, and the (thick) tensor ideal generated by $St_r$ contains exactly the indecomposable tilting modules $T_i$ with $i\ge p^r-1$.

For any $r\ge1$, the \emph{higher Verlinde category} $\Ver_{p^r}$ is defined as the monoidal abelian envelope of the quotient category
$$
\Tilt\SL_2 / \langle St_r \rangle ,
$$
of the category of tilting modules by the tensor ideal generated by $St_r$. The quotient itself is an additive and idempotent complete linear rigid symmetric monoidal category (a \emph{pseudo-tensor category}). For $r=1$, the quotient happens to be semisimple, and hence a tensor category, the Verlinde category $\Ver_p$. For $r>1$, the existence of the abelian envelope is a non-trivial fact established in \cites{BEO,Cou-abelian}. All higher Verlinde categories are tensor categories.

\medskip

By construction, there is a linear symmetric monoidal functor
$$
F_r\: \Tilt\SL_2 \to \Tilt\SL_2/ \langle St_r \rangle \to \Ver_{p^r} ,
$$
the quotient functor followed by the abelian envelope. We will use the following extension of this functor. For all $r\ge1$, define the full subcategory
$$
\ol\cT_r \:= \{ X\in\Rep\SL_2: St_{r-1}\o X\in\Tilt\SL_2\} \quad \subset\Rep\SL_2 .
$$
By construction, $\ol\cT_r$ is an additive idempotent complete rigid symmetric monoidal subcategory of $\Rep \SL_2$, containing $\Tilt \SL_2$. Clearly, $\Tilt\SL_2\subset\ol\cT_r$.

\begin{lemma}[\cite{Newton}*{Section~3}] There is a unique linear symmetric monoidal functor $\ol F_r\:\ol\cT_r\to\Ver_{p^r}$ that agrees with $F_r$ on $\Tilt\SL_2$.
\end{lemma}

By \cite{BEO}*{Section~4.10}, there is also a tensor functor
$$
\Ffrak_r \: \Ver_{p^r} \to \Ver_{p^{r+1}} ,
$$
corresponding, in some sense, to the Frobenius twist in $\Rep\SL_2$. This functor is particularly useful to describe the simple objects in $\Ver_{p^r}$. 

The simple objects $L_i\in\Ver_{p^r}$ can be labeled by integers $0\le i\le p^{r-1}(p-1)-1$. We have $L_0=\one$, the tensor unit, and if $p>2$, we have $L_{p^{r-1}(p-2)}=:\ol\one$, an object tensor-squaring to $\one$. These are the only invertible objects. The simples in $\Ver_{p^r}$ and those in $\Ver_{p^{r+1}}$ are related by
$$
L_a \o \Ffrak_r L_b \cong L_{pb+a}\quad
\text{for all }0\le a\le p-1, 0\le b\le p^{r-1}(p-1)-1.
$$
In particular, all simple objects in $\Ver_{p^{r+1}}$ arise in this way. Note that in the formula, $L_b$ is a simple object in $\Ver_{p^r}$, while $L_a$ and $L_{pb+a}$ are simples in $\Ver_{p^{r+1}}$. In general, we have to specify the parent category for a simple object of the form $L_i$ in order to identify it uniquely. The simples are self-dual, by this formula and the self-duality of the $\SL_2$-tilting modules. Hence, the non-vanishing degrees of $\wedge^d L_i$ and $\Lambda^d L_i$ agree.

\section{Degree bounds in higher Verlinde categories}

\subsection{Transfer functors} \label{sec:transfer} 

We fix $r\ge1$ and $1\le a\le p-1$.  The main goal of this section is to define functors that relate (transfer) symmetric powers in $\Ver_{p^{r+1}}$ to symmetric powers in $\Rep\GL_{a+1}\boxtimes \Ver_{p^r}$. The precise statement is contained in \Cref{prop:G-tensor} below. We start by introducing the main players.

Let
$$ 
\cB_{r,a} \subset \Rep\GL_{a+1}
$$
be the full subcategory of $\GL_{a+1}$-modules whose weights $\mu=(\mu_1,\dots,\mu_{a+1})$ satisfy
$0\le \mu_i \le p^r$  for all $1\le i\le a+1$; this is a finite Serre subcategory. Let
$$
\phi\: \SL_2 \to \GL_{a+1}
$$
be the $\SL_2$-representation $T_a$, i.e., the homomorphism such that $\phi^*V\cong T_a$, where $V\in\Rep\GL_{a+1}$ is the natural $(a+1)$-dimensional representation. Note that $\phi$ is \emph{principal}, see \cite{Newton}*{Section~1.4}, as any matrix $\left(\begin{smallmatrix} 1&0 \\ x&1\end{smallmatrix}\right)$ with $x\neq0$ acts by a Jordan block of size $a+1$ on $T_a=\Sym^a T_1$. Let
$$
\cD_{r,a} \subset \Rep\GL_{a+1}
$$
be the full subcategory of objects $X$ such that $\phi^*X\in\ol\cT_{r+1}$ (with $\ol\cT_{r+1}$ as defined in \Cref{bg:ver}); it is an additive and idempotent complete monoidal subcategory.

\begin{lemma} \label{lem:B-in-D} We have an inclusion $\cB_{r,a}\subset\cD_{r,a}$ of full subcategories of $\Rep \GL_{a+1}$.
\end{lemma}

\begin{proof}
Let $Q:=L((p^r-1)(a,a-1,\dots,0))$ be a choice of a simple $\GL_{a+1}$-representation that restricts to the $r$-th Steinberg representation of $\SL_{a+1}$. We first claim that $Q\o X$ is tilting for all $X\in\cB_{r,a}$.

As tilting representations have trivial extensions, it suffices to show this for $X=L(\lambda)$, a simple $\GL_{a+1}$-representation in $\cB_{r,a}$ with highest weight $\lambda=(\lambda_1,\dots,\lambda_{a+1})$. Let $\alpha_0$ be the highest root of $\SL_{a+1}$, then as $L(\lambda)\in\cB_{r,a}$, 
$$
\langle \lambda,\alpha_0^\vee \rangle = \lambda_1-\lambda_{a+1} \le p^r .
$$
Hence, by \cite{KN}*{Theorem~7.2.2}, $Q\o X$ has a good filtration as an $\SL_{a+1}$-representation. This is a good filtration even as a $\GL_{a+1}$-representation as $Q$ and $X$ are simple, and hence the center acts by a single character on $Q\o X$. The same argument for $X^*$ with highest weight $(-\lambda_{a+1},\dots,-\lambda_1)$ shows that $Q\o X^*$ has a good filtration. As $Q$ is self-dual up to a power of the determinant representation, this implies that $Q\o X$ also has a Weyl filtration. Hence, $Q\o X$ is tilting, establishing the claim.

We now claim that $\phi^*$ preserves tilting representations.  Indeed, indecomposable tilting representations for $\GL_{a+1}$ are direct summands of tensor products involving exterior powers of $V$ of degree at most $a$ and determinant twists. Such exterior powers are direct summands of tensor powers of $V$, as $a<p$, and determinant twists restrict trivially to $\SL_2$. Hence, for the second claim it suffices to note that $\phi^*V\cong T_a$ is a tilting $\SL_2$-representation, by definition.

The two claims together show that $\phi^*Q\o\phi^*X$ is tilting. If $p>2$, the argument in \cite{Newton}*{Lemma~2.12} shows that $\phi^*Q$ is in the thick tensor ideal generated by $T_{p^r-1}$, but not in the one generated by $T_{p^{r+1}-1}$. If $p=2$, then $a=1$, $\phi$ is the inclusion $\SL_2\to\GL_2$, and $\phi^*Q=St_r$, so the same conclusion holds. Hence, by one of the equivalent definitions of $\ol\cT_{r+1}$ in \cite{Newton}*{Section~3.1}, $\phi^*X$ lies in $\ol\cT_{r+1}$.
\end{proof}

Let $\cP_r$ be the linear category of projective objects in $\Ver_{p^r}$. Let $[\cA,\cB]$ denote the category of linear functors between two linear categories $\cA,\cB$. Such categories of functors can be used to realize Deligne tensor products. Concretely, for any locally finite abelian linear category $\cA$, by \cite{EGNO}*{Sections~1.8, 1.11}, we have an equivalence
$$
\cA\boxtimes \Ver_{p^r} \simeq [\cP_r\op,\cA],\quad
X\boxtimes Y \mapsto (P\mapsto \Hom(P,Y)\o X) .
$$
Using such equivalences, the embeddings $\cB_{r,a}\subset\cD_{r,a}\subset\Rep\GL_{a+1}$ can be turned into embeddings
$$
\cB_{r,a}\boxtimes\Ver_{p^r} \subset [\cP_r\op,\cD_{r,a}] \subset \Rep\GL_{a+1}\boxtimes\Ver_{p^r} .
$$

Let 
$\ol F:=\ol F_{r+1}\: \ol\cT_{r+1} \to \Ver_{p^{r+1}}$ and
$\Ffrak:=\Ffrak_r\: \Ver_{p^r} \to \Ver_{p^{r+1}}$
be the functors discussed in \Cref{bg:ver}.

\begin{definition} Define linear functors
$$
G_0\: \cD_{r,a} \to \Ver_{p^{r+1}} ,\quad
X\mapsto \ol F(\phi^*X),
$$
$$
G\: [\cP_r\op,\cD_{r,a}] \xrightarrow{G_0} [\cP_r\op,\Ver_{p^{r+1}}] \simeq \Ver_{p^{r+1}}\boxtimes\Ver_{p^r} \xrightarrow{\o(\Id\boxtimes\Ffrak)} \Ver_{p^{r+1}} ,\quad
X\boxtimes Y \mapsto G_0(X)\o\Ffrak Y ,
$$
for all $X\in\cD_{r,a}$, $Y\in\Ver_{p^r}$.
\end{definition}

Here, the formula describing $G$ on objects is interpreted by realizing $[\cP_r\op,\cD_{r,a}]$ as a subcategory of $\Rep\GL_{a+1}\boxtimes\Ver_{p^r}$ with terms in $\cD_{r,a}$ (or alternatively, we extend the notation $X\boxtimes Y$ for $\Hom(-,Y)\otimes X$ even when the target is not abelian, and so the functor category is not a Deligne tensor product).

The most important properties of these functors can be summarized as follows.

\begin{lemma} \label{lem:G-G0}
(a) $[\cP_r\op,\cD_{r,a}]$ is an additive and idempotent complete symmetric monoidal subcategory of $\Rep\GL_{a+1}\boxtimes\Ver_{p^r}$.

(b) $G_0$ and $G$ are symmetric monoidal.

(c) $G_0$ and $G$ send finite exact sequences in the abelian ambient categories of their sources (i.e., $\Rep\GL_{a+1}$ and $\Rep\GL_{a+1}\boxtimes\Ver_{p^r}$, respectively) with terms in the respective sources to exact sequences in the target. In particular, the restriction of $G$ to $\cB_{r,a}\boxtimes\Ver_{p^r}$ is exact.
\end{lemma}

\begin{proof} 
(a) As the tensor product functor $\Ver_{p^r}\boxtimes\Ver_{p^r}\to\Ver_{p^r}$ is exact, for any projective object $P\in\cP_r$, the functor $\Hom(P, (-\otimes-))$ is represented by a projective object $P'$ in $\Ver_{p^r}\boxtimes\Ver_{p^r}$. Pick $(P_i)_i$ and $(Q_i)_i$ in $\cP_r$ such that
$$
P'\cong \bigoplus_i P_i \boxtimes Q_i
$$
(which is possible as $\kk$ is algebraically closed). 
Then for any $F,F'$ in $[\cP_r\op,\cD_{r,a}]$,
$$
(F\o F')(P)\cong \bigoplus_i F(P_i)\o F'(Q_i) \quad \in \cD_{r,a}
$$
using the discussed closure properties of $\cD_{r,a}$. This shows that $[\cP_r\op,\cD_{r,a}]$ is closed under tensor products. As the tensor unit of $\GL_{a+1}$ lies in $\cD_{r,a}$ and as $\cD_{r,a}$ is closed under direct sums and summands, this shows the assertion.

(b) $G_0$ is symmetric monoidal by definition, which then also makes $G$ symmetric monoidal.

(c) As $\phi^*$ is exact, exactness of $G_0$ in the asserted sense follows from \cite{Newton}*{Proposition~4.8}. This implies exactness of $G$, as Deligne tensor products of exact functors on finite abelian categories are exact, and the tensor product in $\Ver_{p^{r+1}}$ is biexact.
\end{proof}

Recall that $V\in\Rep\GL_{a+1}$ is the natural representation. The main feature of the functors $G$ and $G_0$ relevant for us is that they play well with symmetric powers, in the following sense. 

\begin{proposition} \label{prop:G-tensor} \quad
(a) Assume $X\in\cD_{r,a}$, $d\ge2$, $\Gamma^2X\in\cD_{r,a}$, and $\Sym^i X,\wedge^i X\in\cD_{r,a}$ for all $0\le i\le d$. Then
$$
G_0( \Sym^i(X)) \cong \Sym^i (G_0(X)) 
\quad
\text{for all }0\le i\le d.
$$
Moreover, the condition on $\Gamma^2X$ holds automatically for $p>2$.

(b) Let $Y\in\Ver_{p^r}$ be an object with $m(Y)\le p^r-2$ and let $d\ge2$ be an integer. Then for $k=0,1,2$, the objects 
$$
(V\boxtimes Y)^{\o k}\o \Sym^{d-k}(V\boxtimes Y)
$$
computed in $\Rep\GL_{a+1}\boxtimes\Ver_{p^r}$ lie in $\cB_{r,a}\boxtimes\Ver_{p^r}$. The same holds with $\Sym$ and $m(Y)$ replaced by $\wedge$ and $n(Y^*)$.

(c) For any $Y\in\Ver_{p^r}$ with $m(Y)\le p^r-2$ and $d\ge0$, we have an isomorphism
$$
G(\Sym^d(V\boxtimes Y))\cong \Sym^d(L_a\o \Ffrak Y) ,
$$
compatible with the multiplication $\mu_Z\:Z\o\Sym^d Z\to\Sym^{d+1} Z$ applied on each side for the respective objects $Z$ whose symmetric powers occur. The same holds with $\Sym$ replaced by $\wedge$, assuming $n(Y^*)\le p^r-2$.
\end{proposition}

\begin{proof} For any object $Z$ in a tensor category with symmetric braiding $c$ let $\iota_Z$ be the embedding of $\Lambda^2 Z=\im(\id-c_{Z,Z})$ in $Z\o Z$.
Recall from \cite{CEN}*{Lemma~5.2.3} that, for $d\ge2$, we have a presentation of $\Sym^d Z$ given by (right exactness of) the following sequence:
\begin{equation}
\Lambda^2 Z\o \Sym^{d-2} Z 
 \xrightarrow{(Z\o\mu^Z_{d-1})(\iota_Z\o\id)}
Z \o \Sym^{d-1} Z
 \xrightarrow{\mu^Z_d}
\Sym^d Z
 \to 0.
\label{eq:presentation}
\end{equation}
Similarly, when considering the versions for exterior powers instead of symmetric ones, recall that we have a presentation
\begin{equation}
\Gamma^2 Z\o \wedge^{d-2} Z 
 \xrightarrow{(Z\o\mu^Z_{d-1})(\iota_Z\o\id)}
Z \o \wedge^{d-1} Z
 \xrightarrow{\mu^Z_d}
\wedge^d Z
 \to 0,
\label{eq:presentation-wedge}
\end{equation}
where now $\iota_Z$ is the embedding of $\Gamma^2Z=\ker(\id-c_{Z,Z})$ in $Z\o Z$, and where we denote the two embeddings and multiplication morphisms by the same symbols.

(a) This is the argument as in \cite{Newton}*{Proposition~4.9}: First we note that all terms of the exact sequence
$$
0 \to \Gamma^2 X \to X^{\o 2} \xrightarrow{\id-c_{X,X}} X^{\o 2} \to \Sym^2 X \to 0
$$
lie in $\cD_{r,a}$, so $G_0$ sends it to an exact sequence by \Cref{lem:G-G0}(c). As $\Lambda^2 X\cong \wedge^2 X$, $G_0$ sends the exact sequence
$$
0 \to \Lambda^2 X \to X^{\o 2} \to \Sym^2(X) \to 0
$$
to an exact sequence. As $G_0$ is symmetric monoidal, this shows 
$$
\Sym^2 G_0(X)\cong G_0(X^{\o 2})/G_0(\Lambda^2 X) .
$$  

Now specialize $Z$ to $X\in\Rep\GL_{a+1}$ in \Cref{eq:presentation}. This produces a part of the ordinary Koszul complex, which is a finite exact sequence. All terms lie in $\cD_{r,a}$ by the assumptions. Hence, by \Cref{lem:G-G0}(c), the image of this sequence under $G_0$ remains exact, which implies the assertion.

The condition on $\Gamma^2 X$ is automatic for $p>2$, since it is a direct summand of $X^{\o 2}$ given by a symmetrizer in this case.

(b) As a representation of the torus $T\subset \GL_{a+1}$, $V$ is a direct sum of one-dimensional weight spaces of weights  represented by elementary vectors. Hence, in $\Rep T\boxtimes\Ver_{p^r}$, the object $\Sym^{d-k}(V\boxtimes Y)$ is a direct sum of weight spaces (i.e., weight objects of $\Ver_{p^r}$) of the form 
$$
\bigotimes_{i=1}^{a+1} \Sym^{m_i} Y ,\quad
\text{for } m_1+\dots+m_{a+1} = d-k,
$$
of weight $(m_1,\dots,m_{a+1})$, using \Cref{eq:direct-sums}. This shows that the weights of $\Sym^{d-k}(V\boxtimes Y)$ are componentwise bounded by $m(Y)$, in particular, by $p^r-2$. Similarly, the weights of $\wedge^{d-k}(V\boxtimes Y)$ are componentwise bounded by $n(Y^*)$, and the weights of $(V\boxtimes Y)^{\o k}$ are componentwise at most $k$, in particular, at most $2$.

(c) Specialize $Z$ to $V\boxtimes Y\in \Rep\GL_{a+1}\boxtimes\Ver_{p^r}$ in \Cref{eq:presentation}. Set $W:=L_a\o \Ffrak Y$.

As $p^r\ge 2$, $Z^{\o 2}$ lies in $\cB_{r,a}\boxtimes\Ver_{p^r}$. As the latter is a Serre subcategory of $\Rep\GL_{a+1}\boxtimes\Ver_{p^r}$, $G$ is exact on it by \Cref{lem:G-G0}(c). As $G$ is symmetric monoidal by \Cref{lem:G-G0}(b),
$$
G(\Lambda^2 Z)\cong \Lambda^2 G(Z) 
\quad\text{ and }\quad
G(\Gamma^2 Z)\cong \Gamma^2 G(Z) ,
$$
compatibly with the embeddings into $Z^{\o 2}$ and $G(Z)^{\o 2}$, respectively.

We can now show the assertion by induction in $d$.
The cases $d=0,1$ hold by definition of $G$. Assume by induction the statement is true up to degree $d-1$. Then we can apply the functor $G$ to the presentation of $Z$ as in \Cref{eq:presentation}, where it is exact by \Cref{lem:G-G0}(c), using that by part (b), all relevant objects lie in $\cB_{r,a}\boxtimes\Ver_{p^r}$. Note that the first term is a subobject of $Z^{\o2}\o\Sym^{d-2}Z$. Using the symmetric and monoidal structure and the induction hypothesis, we obtain the corresponding presentation for $\Sym^d W$, i.e., the desired isomorphism at degree $d$, compatible with multiplication by the cokernel construction.

The proof for exterior instead of symmetric powers is similar, using \Cref{eq:presentation-wedge} instead of \Cref{eq:presentation}.
\end{proof}

\subsection{Degree computations} As a first application, the transfer functors discussed in \Cref{sec:transfer} yield a recursion formula for symmetric and exterior degrees.

\begin{proposition} \label{prop:recursion} Assume $Y\in\Ver_{p^r}$ is simple and $1\le a\le p-1$. Then
$$
m(Y)\le p^r-2 
\quad\Rightarrow\quad
m(L_a \o \Ffrak Y) = (a+1) m(Y) 
\quad\text{ and }
$$
$$
n(Y) \le p^r-2 
\quad\Rightarrow\quad
n(L_a \o \Ffrak Y) = (a+1) n(Y)
\quad\text{in } \Ver_{p^{r+1}} .
$$
\end{proposition}

\begin{proof} Assume $m(Y)\le p^r-2$. We set $m:=m(Y)$, $t:=a+1$, and $Z:=V\boxtimes Y$ in $\Rep\GL_t\boxtimes\Ver_{p^r}$. Applying the forgetful functor $\Rep\GL_t\to\Vect$ on the first tensor factor, $Z$ becomes $Y^t$. Using the direct sum formula \Cref{eq:direct-sums}, we compute
$$
W:=\Sym^{tm}(Y^t) \cong (\Sym^m Y)^{\o t} ,
$$
which is invertible, and hence also simple. This means that $\Sym^{tm} Z$ is of the form $\chi\boxtimes W$, for some algebraic character $\chi$ of $\GL_t$. As $\SL_2$ has no non-trivial algebraic character, \Cref{prop:G-tensor}(c) now shows that
$$
\Sym^{tm}(L_a\o \Ffrak Y) \cong G(\Sym^{tm}(Z)) \cong G(\chi\boxtimes W) \cong \Ffrak W
$$
is invertible. $L_a$, and hence $L_a\o \Ffrak Y$, is not invertible. This shows that $m(L_a\o \Ffrak Y)=tm$, as asserted.

The assertion on exterior degrees follows completely analogously, since both $Y$ and $L_a\o\Ffrak Y$ are self-dual.
\end{proof}

As a second application, the transfer functor also enables us to compute some explicit symmetric and exterior degrees in $\Ver_{p^r}$, more specifically for simple objects $L_b$ with small indices $b$.

\begin{proposition} \label{prop:restricted} We have 
$$
n(L_b) = b+1 ,\quad
m(L_b) = p^r - (b+1) \quad
\text{in } \Ver_{p^r} ,
$$
if $r=1$ and $1\le b\le p-3$, or if $r\ge2$ and $1\le b\le p-1$.
\end{proposition}

\begin{proof} For $b=1$, the assertion is \cite{BE-cohomology}*{Beginning of Section~6, Corollary~6.4}, so we may assume $p>2$.

For $r=1$, computing the exterior power of $\SL_2$-representations shows $n(L_b)=b+1$, and $\ol\one\o L_b\cong L_{p-2-b}$ gives $m(L_b)=p-b-1$. Now assume $r\ge2$. For the exterior degrees, using \Cref{prop:recursion} with $Y=\ol\one\in\Ver_{p^{r-1}}$, we get
$$
n(L_b) = m(L_b\o \Ffrak \ol\one) = b+1 ,
$$
as asserted.

For the symmetric degrees, first note that 
$T_b \cong (\Sym^b T_1)^*\cong\Gamma^b T_1$ by self-duality of $T_b$ and $T_1$. Let $V$ be the natural $(b+1)$-dimensional representation of $\GL_{b+1}$, so that $\phi^*V\cong T_b$. By a modular version of Hermite reciprocity \cite{RS-hermite}*{Theorem~2.1}, for $i\geq0$ we have
$$
\phi^*(\Sym^i V)
 \cong \Sym^i(\Gamma^b T_1)
 \cong \Lambda^b(\Sym^{i+b-1} T_1).
$$
Now let $d=p^r-(b+1)$ and suppose $0\leq i\leq d$. This means $\Sym^{i+b-1} T_1$ has weights bounded by $i+b-1\le p^r-2$, so it is in $\ol\cT_r$ by \cite{BEO}*{Lemma~3.3} and vanishing of extensions between tilting modules. Consequently, the right hand side is in $\ol\cT_r$, since the exterior power is given by an antisymmetrizer and $\ol\cT_r$ is idempotent complete. This shows that $\Sym^i V$ lies in $\cD_{r-1,b}$. Moreover, any $\GL_{b+1}$-weight of $\wedge^i V$ has coordinates $0$ or $1$, so $\wedge^i V$ lies in $\cB_{r-1,b}\subset\cD_{r-1,b}$, using \Cref{lem:B-in-D}. So $V$ and $d$ satisfy the conditions of \Cref{prop:G-tensor}(a), and thus
$$
\Sym^d L_b
 \cong \ol F_r(\Sym^d T_b)
 \cong \ol F_r(\Lambda^b(\Sym^{p^r-2} T_1)).
$$
We can then swap $\ol F_r$ and $\Lambda^b(-)$, again using the antisymmetrizer description of $\Lambda^b(-)$. By \cite{BE-cohomology}*{Corollary~6.4(i)} and \cite{Newton}*{Remark~4.10}, we have $\ol F_r(\Sym^{p^r-2}T_1)\cong\ol\one$, so the right-hand side is $\ol\one^{\o b}$ which is invertible. As $L_b$ is not invertible, this shows that $m(L_b)=d$ by \Cref{lem:invertible}, as asserted.
\end{proof}

Taken together, \Cref{prop:recursion} and \Cref{prop:restricted} yield a complete description of the symmetric and exterior degrees of simple objects in all higher Verlinde categories.

\begin{theorem} \label{thm:formulas} Assume $p>2$. Let $L_i$ be a non-invertible simple object of $\Ver_{p^r}$, i.e., $1\le i\le p^{r-1}(p-1)-1$ and $i\neq p^{r-1}(p-2)$. Let $i=i_0+pi_1+\dots+p^{r-1} i_{r-1}$ be the $p$-adic expansion. Set $j:=i+p^{r-1}(p-2-2 i_{r-1})$. Let $k$ be the maximal index such that $i_k>0$. Set $m:=r-1$ if $i_{r-1}\notin\{0,p-2\}$; otherwise, let $0\le m\le r-2$ be maximal such that $i_m>0$. Then
$$
m(L_i) = (p^{r-k}-i_k-1)\prod_{\ell=0}^{k-1} (i_\ell+1) ,\quad
n(L_i)=m(\ol\one\o L_i)=m(L_j)$$
$$
m(L_i)+n(L_i) = p^{r-m} \prod_{\ell=0}^{m-1}(i_\ell+1) \le p^r .
$$
\end{theorem}

\begin{proof} First note that all maxima are well-defined for those $i$ we consider. The assertion on exterior degrees now follows from the one for symmetric degrees, as $\ol\one\o L_i\cong L_j$ by the $\SL_2$ tensor product rule in $\Ver_p$ and the tensor product formula in \Cref{bg:ver}, in particular, $\ol\one\cong L_{p^{r-1}(p-2)}$.

For $r=1$, the assertion on the symmetric degrees and on the sum of the degrees is \Cref{prop:restricted}. Assume all formulas hold in $\Ver_{p^r}$, for some $r\ge1$. Any simple object in $\Ver_{p^{r+1}}$ is of the form $W=L_a\o \Ffrak Y$ for $0\le a\le p-1$ and $Y$ a simple of $\Ver_{p^r}$. If $a=0$, $L_a$ is the tensor unit, so $m(W)=m(Y)$ and $n(W)=n(Y)$. If $Y$ is invertible, then it is the tensor unit or the odd tensor unit, in which case $(m(W),n(W))=(m(L_a),n(L_a))$ or $(n(L_a),m(L_a))$, respectively. In the remaining case, we have $1\le a\le p-1$ and $Y\in\Ver_{p^r}$ not invertible. Hence,
$$
m(Y)+n(Y)\le p^r ,  \quad
m(Y),n(Y)\ge 2 ,
$$
which implies $m(Y),n(Y)\le p^r-2$. Now from \Cref{prop:recursion}, we get
$$
m(W) = (a+1)m(Y),\quad 
n(W) = (a+1)n(Y) .
$$
In all cases, we have determined the symmetric and exterior degree of $W$, and it can be checked that they are given by the asserted formulas.
\end{proof}

\begin{theorem} \label{thm:formula-p2} Assume $p=2$ and $r\ge2$ (recall that $\Ver_2\simeq\Vect$). Let $L_i$ be a non-invertible simple of $\Ver_{p^r}=\Ver_{2^r}$, i.e., $1\le i\le 2^{r-1}-1$. Let $i=i_0+2i_1+\dots+2^{r-1} i_{r-1}$ be the $2$-adic expansion. Let $k$ be the maximal index such that $i_k>0$. Let $s$ be the number of non-zero binary digits of $i$.
Then
$$
m(L_i) = 2^s (2^{r-k-1}-1) ,\quad
n(L_i) = 2^s ,\quad 
m(L_i) + n(L_i) = 2^{r+s-k-1} \le 2^r .
$$
\end{theorem}

\begin{proof} For $i=1$, the assertion follows from \Cref{prop:restricted}. As this is the only choice of $i$ for $r=2$, we have confirmed this case. Assume now the assertion holds for some $r\ge2$, and consider a non-invertible simple $L_i$ in $\Ver_{p^{r+1}}$, i.e.,
$$
L_i \cong L_a \o \Ffrak L_b 
,\quad\text{ for some } a\in\{0,1\}, L_b\in\Ver_{p^r}\text{ simple }.
$$
If $a=0$, then $L_i=\Ffrak L_b$, so the assertion follows from the induction hypothesis. If $L_b$ is invertible, then $L_b\cong\one$ and $L_i\cong L_a$, so \Cref{prop:restricted} applies. Otherwise, $a=1$ and $L_b$ is non-invertible, so using the induction hypothesis
$$
m(L_b)+n(L_b)\le p^r ,\quad m(L_b),n(L_b)\ge 2 ,
$$
which implies that $m(L_b),n(L_b)\le p^r-2$. Hence \Cref{prop:recursion} applies. In any case, the asserted formula follows.
\end{proof}

\begin{corollary} \label{cor:bounds} All simple objects of $\Ver_{p^r}$ not isomorphic to $\one$ have finite symmetric degrees, and all non-invertible simple objects $Y$ satisfy $m(Y)+n(Y)\le p^r$.
\end{corollary}

\begin{remark} \label{rem:correction}
It was conjectured in \cite{CEN}*{Conjecture~5.1.3} that equality holds in \Cref{cor:bounds} for $r\ge2$ and $Y$ outside the image of $\Ver_{p^{r-1}}$, i.e., $m(Y)+n(Y)=p^r$. This is not true, for instance, if $p>2$ and $p^{r-1}>3$, then \Cref{prop:recursion} yields for $L_{p+1}\cong L_1\o \Ffrak L_1$ that
$$
m(L_{p+1}) + n(L_{p+1}) = (2p^{r-1} - 4) + 4 = 2p^{r-1} \neq p^r .
$$
On the other hand, the bound in \Cref{cor:bounds} is sharp and attained by $L_1$ for every $r\ge2$.
\end{remark}

\section{Maximal nilpotence and its consequences}
\label{sec:MN}

Recall from \cite{Cou-commalg}*{Section~3.2} that a tensor category $\cC$ is called \emph{maximally nilpotent} if the following two conditions are satisfied:
\begin{itemize}
\item[\textbf{(MN1)}] For every simple object $L\not\cong\one$ in $\cC$ there is a $d\ge1$ such that $\Sym^d L=0$.
\item[\textbf{(MN2)}] For every non-split mono $\one\hookrightarrow X$ in $\cC$, there is a $d>0$ such that the induced morphism $\one\to\Sym^d X$ is zero.
\end{itemize}

\begin{proposition} \label{cor:MN}
$\Ver_{p^r}$ is maximally nilpotent for all primes $p$ and all $1\le r\le \infty$.
\end{proposition}

\begin{proof} For $p=2$, this is \cite{CEO-incompressible}*{Theorem~9.2.1}, so assume $p$ is odd. Then \Cref{cor:bounds} implies (MN1) for all non-invertible simple objects in $\Ver_{p^r}$. (MN1) is immediate for the odd unit, which is the only non-trivial invertible object. By a combination of \cite{CEO-incompressible}*{Proposition~9.1.1} and \cite{CEO-incompressible}*{Theorem~7.2.2}, (MN1) implies maximal nilpotence for the categories $\Ver_{p^r}$ with $r<\infty$. As $\Ver_{p^\infty}$ is a union of these categories, (MN1) and (MN2) follow for this category, as well.
\end{proof}

We record some of the many consequences of this result. For the rest of the paper, $\kk$ is an algebraically closed field of positive characteristic $p$. All tensor categories will be over $\kk$. Of particular importance will be the case of \emph{finite} tensor categories, i.e., those with finitely many simple objects up to isomorphisms and with enough projectives. By the main result of \cite{CF-beo}, those tensor categories admit a tensor functor to $\Ver_{p^r}$, for some $r\ge1$.

\subsection{Incompressible tensor categories} \label{sec:incompressible}

Recall that a tensor category is called \emph{subterminal} if any two parallel tensor functors to that category are isomorphic.

\begin{proposition} \label{prop:subterminal} $\Ver_{p^r}$ is subterminal for all primes $p$ and $1\le r\le\infty$.
\end{proposition}

\begin{proof}
This is \cite{CEO-incompressible}*{Theorem~7.3.1} using \Cref{cor:MN} and \cite{Cou-gr}*{Theorem~A}.
\end{proof}

Recall that a tensor category $\cC$ is called \emph{incompressible} if any tensor functor $\cC\to\cD$ is an embedding, i.e., an equivalence between $\cC$ and a tensor subcategory of $\cD$. Our results also give an intrinsic characterization of incompressible finite tensor categories.

\begin{proposition} \label{prop:incompressible} A finite tensor category is incompressible if and only if it is maximally nilpotent.
\end{proposition}

\begin{proof}
    If a finite tensor category $\cC$ is maximally nilpotent, it is incompressible, by \cite{CEO-incompressible}*{Theorem~7.2.3(2)}. On the other hand, by \cite{CF-beo}*{Theorem~4.7}, the only finite incompressible tensor categories are the tensor subcategories of $\Ver_{p^r}$, which are maximally nilpotent by \Cref{cor:MN}.
\end{proof}

The following reduction argument for the Benson--Etingof--Ostrik conjecture \cite{BEO}*{Conjecture~1.4} is observed in \cite{CF-beo}*{below Proposition~4.16}. Recall that a tensor category is \emph{finitely generated} if it is generated as a tensor category by one object.

\begin{proposition} To show \cite{BEO}*{Conjecture~1.4}, i.e., that any tensor category of moderate growth admits a tensor functor to $\Ver_{p^\infty}$, it suffices to show this is true for any finitely generated tensor category of moderate growth.
\end{proposition}

\begin{proof} First note that any tensor category of moderate growth is the union of its finitely generated subcategories. By the hypothesis, any such subcategory admits a tensor functor to $\Ver_{p^\infty}$. Since $\Ver_{p^\infty}$ is maximally nilpotent and geometrically reductive, the transfinite gluing argument recalled in \cite{CEO-incompressible}*{Proof of Theorem~7.3.1} (see also \cite{CS-flat}*{Appendix~B.3}) applies, allowing the local tensor functors to be chosen with comparison isomorphisms compatible under restriction and composition. This constructs a tensor functor to $\Ver_{p^\infty}$ as desired.
\end{proof}

\subsection{Commutative algebra} In this section, by algebra we mean commutative algebra in the ind-completion of a tensor category.

Recall from \cites{Cou-commalg,Cou-gr} that a tensor category is said to be \emph{geometrically reductive} if for any non-zero morphism $X\to\one$, there is a $d\ge1$ such that the induced morphism $\Sym^d X\to\one$ is split. Recall from \cite{CEO-incompressible} that an algebra (in the above sense) $A$ of a tensor category over $\kk$ is said to be of \emph{finite type}, if the $\kk$-algebra $\Hom(\one,A)$ is finitely generated and $A$ is finitely generated as a module over $\Hom(\one,A)$.

\begin{proposition} \cite{CEO-incompressible}*{Conjecture~9.3.5} holds: Any finite tensor category $\cC$ is geometrically reductive and any finitely generated algebra $A$ of $\cC$ is of finite type.
\end{proposition}

\begin{proof} This follows from \Cref{cor:MN} and the main result of \cite{CF-beo} by \cite{Cou-gr}*{Theorem~4.1.1}.
\end{proof}

\newcommand\Ann{\operatorname{Ann}}

\begin{proposition} The following hold for any tensor category admitting a tensor functor to $\Ver_{p^\infty}$, in particular, for any finite tensor category.

(a) Hilbert Basis Theorem: Every finitely generated algebra, and every finitely generated module over such an algebra, are noetherian.

(b) Nakayama's lemma: If $A$ is an algebra with an ideal $I$ and a finitely generated module $M$ such that $IM=M$, then $\Ann(M)+I=A$.

(c) For any algebra $A$, epimorphic endomorphisms of finitely generated free $A$-modules (of the form $A\o X$ with $X$ in the tensor category) are isomorphisms.

(d) Every inclusion of commutative Hopf ind-algebras is faithfully flat. 
\end{proposition}

\begin{proof} As all these statements are reflected by tensor functors, it suffices to establish them for $\Ver_{p^\infty}$. This follows from \Cref{cor:MN} and geometric reductivity by \cite{Cou-gr}*{Theorem~A}, using \cite{Cou-commalg}*{Lemma~6.4.4, Theorem~6.1.5, Theorem~6.2.2, Theorem~6.3.1}.
\end{proof}

Recall from \cite{Cou-commalg}*{Definition~4.2.1} that a prime ideal of an algebra $A$ of a tensor category $\cC$ is \emph{visible} if it is the kernel of an algebra morphism $A\to S$, where $S$ is a simple algebra (meaning a simple commutative ind-algebra here).

\begin{proposition} If $\cC$ is finite, then every such prime ideal is visible. 
\end{proposition}

\begin{proof} This follows from \Cref{cor:MN} and the main result of \cite{CF-beo} by \cite{Cou-commalg}*{Theorem~4.2.6}.
\end{proof}

We record some consequences for absolutely flat and exact algebras.

\begin{proposition} Let $\cC$ be a tensor category that admits a tensor functor to $\Ver_{p^\infty}$, for instance, any finite tensor category. Let $A$ be a simple algebra of $\cC$. Then:

(a) $A$ is absolutely flat and artinian, in particular, the tensor product functor $-\o_A-$ on $A$-modules is biexact, and the free $A$-module $A\o X$ has finite length as an $A$-module for all $X\in\cC$.

(b) If $\Hom(\one,A)$ is finite-dimensional, then $A$ is exact: its category of finitely generated modules is a tensor category whose ind-completion is the general category of modules (in the ind-completion of $\cC$). If $\cC$ is finite, then any such $A$ lies in $\cC$.
\end{proposition}

\begin{proof} (a) is \cite{CS-flat}*{Corollary~6.6}, as any such category $\cC$ is \emph{well-fibred} in the sense of \cite{CS-flat} by \Cref{cor:MN} and \cite{Cou-gr}*{Theorem~A}.

(b) By part (a), $A$ is absolutely flat and artinian, together with the assumed hom-finiteness, this means we can apply \cite{CEO-incompressible}*{Proposition~6.1.2}. Using additionally that $A$ is simple, it follows that $\Hom(\one,A)\cong\kk$. The last part is \cite{CEO-incompressible}*{Remark~6.1.3(2)}.
\end{proof}

\subsection{Algebraic geometry}

Maximal nilpotence is also a key ingredient in Coulembier's version of algebraic geometry for tensor categories \cite{Cou-algeo}. For instance, recall that to any algebraic group $G$ in $\Ver_{p^r}$ we can associate an ordinary algebraic group $G_0$ (in $\Vect$). Then the following results apply to all higher Verlinde categories, for which maximal nilpotence is now established by \Cref{cor:MN}.

\begin{proposition}[\cite{CS-homogeneous}*{Theorem~7.5} for $\Ver_{p^r}$] Let $G$ be an algebraic group in $\Ver_{p^r}$, for any prime $p$ and $1\le r\le \infty$. Let $H$ be a closed subgroup. Then

(a) $G/H$ exists as an algebraic separated scheme.

(b) $G_0/H_0\to G/H$ is a closed immersion and a universal homeomorphism.

(c) $G/H$ is affine, quasi-affine, or proper, whenever $G_0/H_0$ has the respective property.
\end{proposition}

\begin{proposition}[\cite{CS-homogeneous}*{Corollary~7.10} for $\Ver_{p^r}$] Every normal subgroup of an affine group scheme in $\Ver_{p^r}$ is the kernel of a group homomorphism.
\end{proposition}

\subsection{Symmetric growth}

Recall from \cite{CEO-incompressible} that the \emph{symmetric growth dimension} of an object $X$ in a tensor category is defined as
$$
\sgd(X) := \limsup_{d\to\infty} \frac{\log(\len(\Sym^{\le d}(X)))}{\log(d)} ,
$$
where $\Sym^{\le d}(X):=\bigoplus_{i=0}^d \Sym^i(X)$ and $\len$ denotes the length of an object in an abelian category. The following answers \cite{CEO-incompressible}*{Question~7.1.4(1--3)} for finite tensor categories.

\begin{proposition} Let $X$ be an object of a finite tensor category $\cC$ with a tensor functor $F\:\cC\to\Ver_{p^r}$, for some $r\ge1$. Then
$$
\sgd(X) = \sgd(F(X)) = \dim\Hom(F(X),\one) \le \FPdim(X),
$$
and $\limsup$ can be replaced by $\lim$ in the definition of $\sgd$.
\end{proposition}

\begin{proof}
    This follows from \Cref{cor:MN} by \cite{CF-beo}*{Proposition~4.15}.
\end{proof}

In \cite{CSZ}, it is shown that any finite tensor category $\cC$ has a universal maximal simple commutative algebra object $\cF(\cC)$.

\begin{proposition} \cite{CSZ}*{Conjecture~B.5} holds: For every $X\in\cC$, $\dim\Hom(X,\cF(\cC))=\sgd(X)$.
\end{proposition}

\begin{proof}
As explained in \cite{CSZ}, scalar extension gives a tensor functor
$$
\cC\to \modfg(\cF(\cC)) =: \cC'
$$
valued in the category of finitely generated $\cF(\cC)$-modules, which is an incompressible finite tensor category.  By the main result of \cite{CF-beo}, $\cC'$ embeds into $\Ver_{p^r}$ for some $r\ge1$, so by \Cref{cor:MN}, $\cC'$ is maximally nilpotent. This implies the assertion, as explained in \cite{CSZ}*{right before Conjecture~B.5}. 
\end{proof}

\section{Alternating powers} \label{sec:alternating}

Still $\kk$ is an algebraically closed field of positive characteristic.

We discuss the related notion of alternating powers. Recall from \cite{CEO-frob}*{Section~2.3}, that for any object $X$ of a tensor category over $\kk$ and for any $d\ge0$, the $d$-th \emph{alternating power} is defined as
$$
A^d(X) := \im( a_d ) ,\quad \text{with }
a_d := \sum_{\sigma\in S_d} (-1)^\sigma \sigma \quad \in\End(X^{\o d}) ,
$$
where we use the algebra morphism $\kk S_d\to\End(X^{\o d})$ defined by the symmetric braiding.
Note that we have morphisms
\begin{equation}
\wedge^d X\twoheadrightarrow A^dX\hookrightarrow \Lambda^d X
\label{eq:ext-powers}    
\end{equation}
for any $X,d$ as above.

\begin{proposition} \label{prop:alternating} Let $\cC$ be any symmetric tensor category of moderate growth over $\kk$. Then for any non-zero $X\in\cC$, there is a $d\ge1$ such that $Z:=A^d(X)$ is invertible and $A^i(X)=0$ for all $i>d$. Moreover, $c_{Z,Z}=\id$ and $Z$ has categorical dimension $1$.
\end{proposition}

\begin{remark}
    In \Cref{prop:alternating}, contrary to the corresponding results for $\wedge$, we do not need to exclude invertible $X$. For instance, for $p>2$ and the odd line $\ol{\one}$ in $\SVect$, we have $A^p\ol{\one}=0\not=\wedge^p \ol{\one}$.
\end{remark}

We prepare the proof with the following auxiliary result.

\begin{lemma} \label{lem:Hopf-algebra} Let $H=\bigoplus_{i=0}^N H_i$ be a graded Hopf algebra with invertible antipode in $\cC$ such that $H_0\cong\one$ and $Z:=H_N\neq0$. Then $Z$ is invertible as an object in $\cC$. Moreover, there is a scalar $s\in\{\pm 1\}$ such that $c_{Z,Z}=s \id_{Z\o Z}$ and $S|_Z=s \id_Z$, where $S$ is the antipode of $H$.
\end{lemma}

\begin{proof} This is a slight abstraction of the proof in \cite{CEO-incompressible}*{Lemma~7.4.5}: As $H$ is a finite Hopf algebra with invertible antipode of $\cC$ (i.e., not just in the ind-completion), it is a Frobenius algebra by \cite{Tak}*{Corollary~4.5}. This means there is an invertible object $L\in\cC$ and an isomorphism 
$$
\phi\: H\cong H^* \o L
$$
as left $H$-modules. Using the definition of the module structure on $H^*$, we see that the components of $\phi$ are of the form
$$
\phi_i\: H_i \to \bigoplus_{0\le j\le N-i} H_j^* \o L \quad 
\text{for } 0\le i\le N.
$$
As $\phi$ is an isomorphism, all $\phi_i$ must be monomorphisms. As $L$ is invertible, and hence, simple, $\phi_N$ must be an isomorphism. This shows that $Z=H_N\cong L$ is invertible. In particular, $Z$ and $Z\o Z$ are simple, so there are $s,t\in\kk^\times$ such that $S|_Z=s \id_Z$ and $c_{Z,Z}=t \id_{Z\o Z}$. As the braiding is symmetric, $t\in\{\pm1\}$.

It remains to show $s=t$. We give a direct proof similar to \cite{BKLT}*{Proof~of~ Proposition~4.10}: let $(\eta,\mu,\varepsilon,\Delta)$ be the remaining structural morphisms of $H$, i.e., unit, multiplication, counit, and comultiplication. Let $\mu_{H\o H}:=(\mu\o\mu)(H\o c\o H)$ be the multiplication in $H\o H$. Then using that $\Delta$ is an algebra morphism, i.e.,
$$
\Delta\mu = \mu_{H\o H}(\Delta\o\Delta) ,
$$
we compute
$$
(\mu\o H) (S\o (\Delta\mu)) (\Delta\o H)
= (H\o \mu) (c\o H) (H\o \Delta) .
$$
Let $\iota\:Z\to H$ be the direct sum inclusion. Pre-composing with $\iota\o\iota$ and using that the bialgebra structure respects degrees, we obtain
$$
(S\iota)\o\iota = c_{H,H}(\iota\o\iota) = (\iota\o\iota) c_{Z,Z} , 
$$
where naturality of the braiding gives the last step. This confirms $s=t$.
\end{proof}

\begin{proof}[Proof of \Cref{prop:alternating}] Set $\cC':=\SVect\boxtimes\cC$, a tensor category that can be identified with the category of $\ZZ/2\ZZ$-graded objects of $\cC$ with braiding $c'_{V,W}:=(-1)^{|V||W|} c_{V,W}$ for homogeneous objects $V,W$.
Let $X':=\ol\one\boxtimes X\in\cC'$ be $X$ placed in odd degree.
Consider the pair $(X',c':=c'_{X',X'})$, where $c':=-c_{X,X}$ is the negative of the symmetric braiding on $X$. This is a symmetric braided object of $\cC'$, i.e., $c'$ satisfies the Yang--Baxter equation and symmetry, as $c$ does. Let $H:=B(X',c')$ be the Nichols algebra associated to this braided object, i.e., the image of the graded Hopf algebra morphism
$$
T(X')\to T\sh(X') , \quad
x_1\o\dots\o x_d \mapsto \sum_{\sigma\in S_d} \psi'(\sigma)(x_1\o\dots\o x_d) ,
$$
using the algebra maps $\psi'\:\kk S_d\to\End(X'^{\o d})$ defined by the symmetric braiding $c'$ of $X'$. Here, $T\sh(X')$ denotes the shuffle Hopf algebra. Let $\psi\:\kk S_d\to\End(X'^{\o d})$ be the algebra morphism defined by $c$. Then
$$
\psi'(\sigma) = (-1)^\sigma \psi(\sigma) \quad
\text{for all }\sigma\in S_d.
$$
This means that the graded subspaces $(H_d)_{d\ge0}$ of $H$ are given by the alternating powers, i.e., $H_d\cong \ol\one^{\o d}\boxtimes A^d(X)$ for all $d\ge0$.

By \cite{CEO-frob}*{Corollary~4.11}, there is an $N\ge1$ such that $A^N(X)\neq0$ and $A^i(X)=0$ for all $i>N$. Hence, by \Cref{lem:Hopf-algebra}, $Z':=H_N$ is invertible, so $Z:=A^N(X)$ is invertible, and there is a scalar $s$ such that $c'_{Z',Z'}=s \id_{Z'\o Z'}$ and $S|_{Z'}=s \id_{Z'}$.

Since $H$ is commutative in $\cC'$ and generated by the primitive object $X'$, its antipode acts on $H_N$ by $(-1)^N$. Hence, using $Z'\cong\ol\one^{\o N}\boxtimes Z$ and the definition of $c'$, and identifying $Z'\o Z'$ with $\one\boxtimes(Z\o Z)$,
$$
(-1)^N \id_{Z\o Z} = (-1)^{N^2} c_{Z,Z} ,
$$
which implies $c_{Z,Z}=\id_{Z\o Z}$. This in turn shows that $Z$ is of dimension~$1$.
\end{proof}

\begin{remark}
    (1) The proof above applies to the exterior algebra $\wedge^\bullet(X)$ (this is a {\em pre-Nichols} algebra associated 
    with $X':=\ol\one\boxtimes X$). Thus if $n\in \ZZ_{\ge0}$ is largest
    with $\wedge^n(X)\ne 0$, then $\wedge^n(X)$ is invertible. 
    This result was proved earlier, see \cite{CEN}*{2.3.2}.
    
    (2) Another variation is to replace anti-symmetrizer $a_d$
    with symmetrizer 
    $$
    s_d := \sum_{\sigma\in S_d}  \sigma \quad \in\End(X^{\o d}).
    $$
\end{remark}

\begin{remark} While the alternating degree of an object is bounded by the exterior degree and the dual exterior degree by \Cref{eq:ext-powers}, they do not have to coincide in general. For instance, one can compute
$$
A^7(L_4)=0\not\cong \one\cong \wedge^7(L_4) 
\quad \text{in } \Ver_9.
$$    
\end{remark}

\section{Hilbert series and \texorpdfstring{$p$}{p}-adic dimensions} \label{sec:dim}

\begin{definition} For $X\in\Ver_{p^r}$, set $v(X):=\dim_\kk\Hom(X,\one)\in\ZZ_{\ge0}$. Pick a short exact sequence
$$
0\to X' \to X \to \one^{\oplus v(X)} \to 0 ,
$$
and set $R(X):=\im(\Sym(X')\to\Sym(X))$, a (finite) object in $\Ver_{p^r}$ by \cite{CEO-incompressible}*{Theorem~7.2.2(4)} and \Cref{cor:MN}. Define 
$$
r(X):=\max\{d\ge0: R(X)_d\neq0\}  \quad \in\ZZ_{\ge0}.
$$
\end{definition}

The notation suggests that $v(X)$ `measures' the `vector' part of $X$ and $r(X)$ the `rest'.

 We also write $x_j=[R(X)_j]\in K_0(\Ver_{p^r})$. By \cite{CEO-incompressible}*{Theorem~3.3.2}, we have an equality
\begin{equation}
    \label{eq:origin}
    \sum_{i\ge 0} [\Sym^i X]\, t^i\;=\; (1-t)^{-v(X)}\left(\sum_{j=0}^{r(X)}x_jt^j\right)\quad\mbox{in}\quad K_0(\Ver_{p^r})[[t]].
\end{equation}
In particular, for any group homomorphism $\psi:K_0(\Ver_{p^r})\to R$ to a commutative ring $R$, it follows that the formal power series $\sum_{i\ge 0} \psi([\Sym^i X]) t^i$ is a rational function, with only poles at $t=1$. We can apply this in a number of situations, leading to the resolution of some open questions.

We start with the $p$-adic dimensions $\Dim_{\pm}(X)\in\mathbb{Z}_p$ defined for objects $X$ in symmetric tensor categories $\cC$ in \cite{EHO}. Concretely for any $X\in\cC$, the $p$-adic dimensions $\Dim_+(X),\Dim_-(X)\in\ZZ_p$ are defined by the identities in $\mathbb{F}_p[[t]]$:
$$
\sum_{k\ge0} \dim(\Sym^k X) t^k = (1-t)^{-\Dim_+(X)} ,\quad
\sum_{k\ge0} \dim(\wedge^k X) t^k = (1+t)^{\Dim_-(X)} ,\quad
$$
where $\dim$ denotes the usual categorical dimension and where $(1+x)^a$ is interpreted as $\prod_{i\ge 0} (1+x^{p^i})^{a_i}$ for any $a\in\ZZ_p$ with $p$-adic expansion $a=a_0+p a_1+\cdots$. If $p>2$, then $\Dim_-(X)$ can be identified with $-\Dim_+(\ol{\one}\otimes X)$ (after potentially adjoining an odd tensor unit to $\cC$, i.e., passing to $\SVect\boxtimes\cC$), so we mainly focus on $\Dim_+$.

By the results in the current paper, an affirmative answer to \cite{EHO}*{Question~3.15} would now follow from an affirmative answer to \cite{BEO}*{Conjecture~1.4}. In particular we can resolve the question for finite tensor categories.
\begin{proposition}\label{prop:padic}
\begin{enumerate}
    \item If $\cC$ is a finite tensor category then $\Dim_{\pm}(X)\in\ZZ\subset\mathbb{Z}_p$ for every $X\in\cC$.
    \item If $X\in\Ver_{p^r}$ for some $r\in\ZZ_{>0}$, then
    \[\Dim_+(X)=v(X)-r(X).\]
\end{enumerate}
    
\end{proposition}
\begin{proof}
For part (2), from \Cref{eq:origin}, we obtain
\[(1-t)^{v(X)-\Dim_+(X)}\;=\; \sum_{j=0}^{r(X)}d_jt^j,\]
    where $d_j=\dim R(X)_j$. Since $R(X)$ is a finite graded Hopf algebra with invertible antipode, \Cref{lem:Hopf-algebra} applies. In particular, $d_{r(X)}=\pm 1\not=0$. So both sides are polynomials of degree $r(X)$, which forces $v(X)-\Dim_+(X)=r(X)$.\footnote{This actually tells us that the invertible $R(X)_{r(X)}$ has same parity as $r(X)$ if $p>2$.}

        Since $p$-adic dimensions are preserved by tensor functors, by \cite{CF-beo}*{Theorem~B} we only need to prove part (1) for $\Ver_{p^r}$. So for $p>2$, we are already done, using that $\Dim_-(X)=-\Dim_+(\ol\one\o X)$. For $p=2$,  recall that \cite{CEN}*{Theorem~2.2.5} states that, with $X\in\Ver_{p^r}$ and $X'\subset X$ as above inducing a filtration on the exterior powers, $\gr \wedge X\cong R'\o S'$, where $R'$ is a quotient of $R(X)$ and $S'$ is a quotient of the symmetric algebra $\Sym(\one^{\oplus v(X)})$ by a homogeneous ideal in the subcategory $\Vect$. This implies
$$
(1+t)^{\Dim_-(X)} = H_{R'}(t) (1-t)^{-v(X)} Q(t) ,
$$
where $H_{R'}$ is the Hilbert polynomial of the finite object $R'$ and $Q$ is a polynomial by the Hilbert-Serre theorem. Multiplying both sides by $(1+t)^{v(X)}=(1-t)^{v(X)}$, we see that $\Dim_-(X)\in\ZZ$.
\end{proof}

\begin{example}
    If $m(X)<\infty$, then $\Dim_+(X)=-m(X)$ and if $n(X^\ast)<\infty$, then $\Dim_-(X)=n(X^\ast)$, in particular, for non-invertible simple $L\in\Ver_{p^r}$, we have
    \[\Dim_+(L)=-m(L)\in \mathbb{Z}_{<-1}\quad\mbox{and}\quad \Dim_-(L)=n(L)\in\mathbb{Z}_{>1}.\]
\end{example}

\begin{remark}
    The alternating version of $p$-adic dimension $\Dim_a(X)\in \ZZ_p$ was defined
    in \cite{CEO-frob}*{8.27} via
$$
\sum_{k\ge0} \dim(A^k(X)) t^k = (1+t)^{\Dim_a(X)}.
$$
It follows immediately from Proposition \ref{prop:alternating}
that 
$$
\Dim_a(X)=\max\{d\ge0: A^d(X)\neq0\}.
$$
In particular, $\Dim_a(X)\in \ZZ_{\ge0}$ for any object of a tensor category of moderate growth.
\end{remark}

We can also resolve \cite{EO-Frob}*{Conjecture~A.4} for finite tensor categories. Following \cite{EO-Frob}*{Appendix~A}, for an object $X$ in a finite tensor category $\cC$, set
\[f_X(t)\;=\;\sum_{i\ge 0} \FPdim(\Sym^i X) \,t^i\;\in\;\mathbb{R}[[t]].\]

\begin{proposition}For every finite tensor category $\cC$ and $X\in\cC$, the formal power series $f_X(t)$ is a rational function with only poles at $t=1$.
\end{proposition}
\begin{proof}
Since tensor functors preserve Frobenius-Perron dimension, it suffices to deal with $\cC=\Ver_{p^r}$ by \cite{CF-beo}*{Theorem~B}.
    The latter follows by applying the ring homomorphism $\FPdim: K_0(\Ver_{p^r})\to\mathbb{R}$ to \Cref{eq:origin}.
\end{proof}

As a final special case of \Cref{eq:origin} we mention:
\begin{lemma}
    For any $X,P\in\Ver_{p^r}$ with $P$ projective, the formal power series \[\sum_{i\ge 0} \dim_\kk\Hom(P,\Sym^i X) t^i\in\mathbb{Z}[[t]]\] is a rational function with only poles at $t=1$. 
\end{lemma}

\begin{remark}
    Results like \Cref{prop:padic}(2) apply to all maximally nilpotent tensor categories. However, we do not expect other maximally nilpotent tensor categories to exist besides $\Ver_{p^{\infty}}$ and its tensor subcategories, so we do not formulate statements in that generality.
\end{remark}


\bibliography{bib}{}
\bibliographystyle{plain}

\end{document}